\documentclass{amsart}
\usepackage[margin=1in]{geometry}

\usepackage%
[bookmarksdepth=subsubsection, unicode]%
{hyperref}

\usepackage{bm,booktabs,comment,mathtools,multicol,tikz-cd,amssymb,stmaryrd}

\usepackage[final,notref,notcite]{showkeys}

\newtheorem{result}{Result}[section]
	\newtheorem{corollary}[result]{Corollary}

	\newtheorem{proposition}[result]{Proposition}
	
	\newtheorem{theorem}[result]{Theorem}

\theoremstyle{definition}
	
	\newtheorem{definition}[result]{Definition}

\theoremstyle{remark}
	
	\newtheorem{remark}[result]{Remark}

\DeclareMathOperator\Hom{Hom}

\DeclareMathOperator\Res{Res}

\newcommand\Spec{\mathrm{Spec}}

\newcommand\bbZ{\mathbb{Z}}
\newcommand\cts{\textnormal{cts}}
\newcommand{\ov}{\overline}

\newcommand\acts\curvearrowright

\newcommand\tn{\textnormal}

\title{Generalized Tate--Nakayama duality for complexes of tori over local function fields} 

\author{Peter Dillery}
\address{Department of Mathematics, Oregon State University, 2000 SW Campus Way, Corvallis, Oregon, 97331, USA.}

\begin{document}

\begin{abstract}
We define the cohomology set $H^{1}(\mathcal{E}, Z \to T \to U)$ for a complex $Z \to T \to U$ with $T$, $U$ tori over a local function field $F$ and $Z$ a finite subgroup of $T$, where $\mathcal{E}$ is the gerbe from \cite{Dillery1}. We follow the method of \cite{Kal26b} and prove that the Tate--Nakayama isomorphism \cite[Theorem 4.10]{Dillery1} for $H^{1}(\mathcal{E}, Z \to T)$ extends to $H^{1}(\mathcal{E}, Z \to T \to U)$, which lets one generalize Kaletha's parametrization \cite{Kal22} of supercuspidal $L$-packets to local function fields.
\end{abstract}

\maketitle

\section{Introduction}
\subsection{Motivation and main results}
For a $p$-adic field $F$, Kaletha defined \cite{kaletha16} the Galois gerbe $W$ banded by a profinite group scheme $u$ and predicted \cite[\S 5.4]{kaletha16} that, for a fixed quasi-split connected reductive group $G$ over $F$ with Whittaker datum $\mathfrak{w}$ and fixed finite central subgroup $Z$, one has for each tempered $L$-parameter $\phi \colon W_{F} \times \mathrm{SL}_{2}(\mathbb{C}) \to {^{L}G}$ a Whittaker-normalized bijection
\begin{equation}\label{eq:rigLLC}
    \Pi_{\phi}:=\bigsqcup_{[z] \in H^{1}(u \to W, Z \to G)} \Pi_{\phi}(G_{[z]}) \xrightarrow{\iota_{\mathfrak{w}}} \text{Irr}(\pi_{0}(Z_{\widehat{G}}(\phi) \times_{\widehat{G}} \widehat{G/Z}))
\end{equation}
satisfying the endoscopic character identities. In the above, $H^{1}(u \to W, Z \to G)$ denotes all equivalence classes of $G(\ov{F})$-valued $1$-cocycles of $W$ whose restriction to $u(\ov{F})$ is an algebraic homomorphism $u(\ov{F}) \to Z(\ov{F})$, $G_{[z]}$ is the inner form obtained by taking the image of $[z]$ in $H^{1}(F,G/Z)$, and $\Pi_{\phi}(G_{[z]})$ denotes the fiber of the local Langlands correspondence for $G_{[z]}$ over the $\widehat{G}$-conjugacy class of $\phi$.

Several cases of the above conjecture (usually without the endoscopic character identities) have been proved; a notable recent example is \cite{Kal22}, which proves that each $\Pi_{\phi}$ can be parametrized by the set $\text{Irr}(\pi_{0}(Z_{\widehat{G}}(\phi) \times_{\widehat{G}} \widehat{G/Z}))$ when $\phi$ is a discrete $L$-parameter $W_{F} \to {^{L}G}$, $G$ splits over a tamely ramified extension of $F$, and the residual characteristic $p$ of $F$ does not divide the order of the absolute Weyl group of $G$---this also includes constructing the Langlands correspondence (yielding $[\phi]$) for all representations $\pi$ lying in such $L$-packets, which Kaletha calls non-singular. Given such a parameter $\phi$, one obtains an elliptic maximal torus $S$ of $G$, and a key step (cf. \cite[\S 4.5]{Kal22}) in the $L$-packet parametrization is carefully comparing the irreducible representations of two extensions of a subgroup of $[N_{G}(S)/S](F)$ by $\mathbb{C}^{\times}$, one arising from $\phi$ and the other from a non-singular representation in its $L$-packet. This comparison requires \cite[Proposition 4.5.1, Lemma 4.5.8]{Kal22}, which use the cohomology set $H^{1}(u \to W, Z \to T \to U)$ for a complex $Z \to T \to U$, where $T$ and $U$ are $F$-rational tori. A critical tool for using this cohomology set is an extension of the Tate--Nakayama duality isomorphism \cite[Theorem 4.8]{kaletha16} from $H^{1}(u \to W, Z \to T)$ to $H^{1}(u \to W, Z \to T \to U)$, which is \cite[Proposition 5.7]{Kal26b}.

When $F$ is a local function field the group scheme $u$ is not smooth and to represent the canonical class in $H^{2}_{\text{fppf}}(F,u)$ one has to use an abstract gerbe $\mathcal{E} \to \Spec(F)$ and consider the set $H^{1}(\mathcal{E}, Z \to G)$ of isomorphism classes of $G$-torsors on $\mathcal{E}$ whose natural $u$-action is induced by a homomorphism $u \to Z$. It is shown in \cite{Dillery1} that if one replaces the cohomology set $H^{1}(u \to W, Z \to G)$ with $H^{1}(\mathcal{E}, Z \to G)$ then the bijections \eqref{eq:rigLLC} (and endoscopic character identities) can be formulated over $F$. Moreover, the set $H^{1}(\mathcal{E}, Z \to G)$ can be understood using linear algebra via a generalization \cite[Theorem 4.10, Theorem 5.10]{Dillery1} of the Tate--Nakayama duality isomorphism.

The goal of this paper is to define the cohomology set  $H^{1}(\mathcal{E}, Z \to T \to U)$ and prove:

\begin{theorem}[Theorem \ref{thm:TN}]\label{thm:intro} The Tate--Nakayama duality isomorphism from \cite[Theorem 4.10]{Dillery1} extends to the cohomology set $H^{1}(\mathcal{E}, Z \to T \to U)$.
\end{theorem}

The main consequence is Theorem \ref{thm:scuspidal}, which says that, under the same assumptions as \cite{Kal22}, the construction of the $L$-packet parametrization from \cite{Kal22} extends to the setting of local function fields.

We emphasize that the proof of Theorem \ref{thm:intro} is a direct adaptation of the arguments of \cite[\S 5]{Kal26b} to the setting of \cite{Dillery1}. In practice, this adaptation consists of rewriting Galois cocycle arguments in terms of \v{C}ech cohomology and torsors. We deliberately follow \cite[\S 5]{Kal26b} as closely as possible (highlighting the pieces which are new in our setting) for ease of readability by those already familiar with the $p$-adic setting.

\subsection{Notation and conventions}
Denote by $F$ a local function field of characteristic $p>0$ with fixed algebraic closure $\ov{F}$ containing a separable closure $F^{s}$ and set $\Gamma:=\text{Gal}(F^{s}/F)$.

For an affine group scheme $A$ over $\Spec(F)$ and an $F$-algebra $R$, we set $A(R):=A(\Spec(R))$ and $A_{R}:= A \times_{\Spec(F)} \Spec(R)$. For a homomorphism of $F$-algebras $f \colon R \to R'$ we also use $f^{\sharp}$ to denote the induced morphism $A(R) \to A(R')$. We denote the homomorphism $\ov{F}^{\bigotimes_{F}n} \to \ov{F}^{\bigotimes_{F}n+1}$ mapping into all but the $j$-factor by $p_{1,\dots,\hat{j},\dots,n+1}$. More generally, when $A$ is commutative, one can use these tensor inclusion maps to define a complex
\begin{equation*}
    A(R) \xrightarrow{d} A(R \otimes_{F} R) \xrightarrow{d} A(R \otimes_{F} R \otimes_{F} R) \xrightarrow{d} \dots
\end{equation*}
and denote the induced groups of cocycles, resp. cohomology classes, resp. cochains by $\check{Z}^{i}(R/F, A)$, resp. $\check{H}^{i}(R/F,A)$, resp. $\check{C}^{i}(R/F,A)$. When $A$ is not necessarily commutative, these objects still make sense for $i=0,1$ as pointed sets, and are denoted in the same way. For concreteness, our $d$ are defined by sending $a \in A(R^{\bigotimes_{F} i})$  to $\prod_{1 \leq j \leq i+1} p^{\sharp}_{1,\dots,\hat{j},\dots,i+1}(a)^{(-1)^{j+1}}$.

A gerbe over $\Spec(F)$ is a stack $\mathcal{E} \xrightarrow{\pi} \text{Schemes/$\Spec(F)$}_{\text{fpqc}}$ fibered in groupoids such that any two objects in a fiber $\mathcal{E}(S)$ for an $F$-scheme $S$ are fpqc-locally isomorphic and such that $\pi$ has a section fpqc-locally. We assume in addition that, for every object $X$ in $\mathcal{E}(R)$ for an $F$-algebra $R$, the group $\text{Aut}_{\Spec(R)}(X)$ is commutative. In this case, this latter automorphism functor defines a sheaf of commutative groups over $\Spec(F)$ called the band of $\mathcal{E}$, see e.g. \cite[Lemma 2.27]{Dillery1}---all gerbes considered in this paper are assumed to have commutative band. If $A$ is a commutative group scheme over $F$, we call a pair $(\mathcal{E},\theta)$ an $A$-gerbe over $F$ if $\mathcal{E}$ is a gerbe over $\Spec(F)$ and $\theta$ is an isomorphism from the band of $\mathcal{E}$ to $A$---in all cases considered in this paper, this isomorphism will be fixed throughout and we omit it from the datum. Because $\mathcal{E} \to \Spec(F)$ is fibered in groupoids, it has finite fibered products and we can therefore transfer the fpqc topology from $\text{Schemes/$\Spec(F)$}_{\text{fpqc}}$ to $\mathcal{E}$ as in \cite[Definition 2.22]{Dillery1}. For an affine group scheme $G$ over $F$, we denote by $Z^{1}(\mathcal{E},G)$ the groupoid of $G \times_{\Spec(F)} \mathcal{E}$-torsors on $\mathcal{E}$, which we frequently refer to as simply $G$-torsors. 

For a ring $R$ and an $R$-group scheme $A$, we use $H^{i}_{\text{fppf}}(R,A)$ or just $H^{i}(R,A)$ to denote $H^{i}_{\text{fppf}}(\Spec(R),A)$. When $G$ is a smooth affine group scheme over $F$, the natural map $\check{H}^{i}(F^{s}/F,G) \to H^{i}_{\text{fppf}}(F,G)$ is an isomorphism and the left-hand set can be canonically identified with $H^{i}(\Gamma, G(F^{s}))$ (where $G$ is always assumed commutative if $i\geq 2$).

Throughout the paper there will often be a fixed finite multiplicative group scheme $Z$ defined over $F$, and for any group $G$ defined over $F$ with a central monomorphism $Z \to G$, we set $\bar{G}=G/Z$, and for $g \in G(R)$ we denote by $\bar{g}$ its image in $\bar{G}(R)$.

\subsection{Acknowledgments}
The author thanks Noriyuki Abe and Kazuma Ohara for suggesting that this paper be written.

AI, specifically Claude Fable 5.2, was used only during the final proofreading stages of this paper to help find potential typos and mathematical errors. 

\section{Setup}\label{sec:setup}

\subsection{Gerbes associated to \v{C}ech cocycles}\label{sec:gerbeb}
For a fixed \v{C}ech $2$-cocycle $\xi \in \check{Z}^{2}(\ov{F}/F,A)$ where $A$ is a commutative $F$-group scheme with $H^{1}(\ov{F}^{\bigotimes_{F} n},A)=0$ for all $n\geq 1$, we can define a specific gerbe $\mathcal{E}_{\xi}$ whose fiber over $\Spec(R)$ for an $F$-algebra $R$ is the groupoid whose objects are $\xi$-twisted $A$-torsors, which are defined as pairs $(T,\psi)$, where $T$ is an $A$-torsor over $\Spec(R \otimes_{F} \ov{F})$ and $\psi$ is a ``twisted gluing map,'' which is an isomorphism of $A$-torsors $p_{2}^{*}T \to p_{1}^{*}T$ such that $(p_{12}^{*}\psi) \circ (p_{23}^{*}\psi) \circ (p_{13}^{*}\psi)^{-1} =m_{\xi}$, where $m_{\xi}$ denotes the ``multiplication by $\xi$'' automorphism of $q_{1}^{*}T$, where $q_{1} \colon R \otimes_{F}\ov{F}  \to R \otimes_{F}\ov{F}^{\bigotimes_{F} 3}$ includes $\ov{F}$ into the first $\ov{F}$-factor. Isomorphisms in this groupoid are isomorphisms of torsors over $\Spec(R \otimes_{F} \ov{F})$ compatible with the twisted gluing maps. This extends functorially to the category of schemes over $\Spec(F)$ and defines an $A$-gerbe over $F$ split over $\ov{F}$ representing the class $[\xi] \in \check{H}^{2}(\ov{F}/F, A)$ (this is proved in \cite[Proposition 2.36]{Dillery1}).

The gerbes $\mathcal{E}_{\xi}$ are useful to work with because, for an affine group scheme $G$ over $F$, the groupoid of $G$-torsors on $\mathcal{E}_{\xi}$ may, by \cite[Proposition 2.50]{Dillery1}, be canonically identified with the groupoid consisting of $\xi$-twisted $G$-torsors, which is the groupoid consisting of pairs $(S,\varphi)$, where $S$ is an $A$-equivariant $G$-torsor over $\ov{F}$ and $\varphi$ is an $A$-equivariant ``twisted gluing'' isomorphism $p_{2}^{*}S \to p_{1}^{*}S$ such that $(p_{12}^{*}\varphi) \circ (p_{23}^{*}\varphi) \circ (p_{13}^{*}\varphi)^{-1} =m_{\xi}$ as an automorphism of $q_{1}^{*}S$. Isomorphisms in this category are isomorphisms of torsors compatible with the twisted gluing maps. 

Moreover, if $f \colon A \to B$ is a morphism of group schemes over $F$ and we have $a \in \check{Z}^{2}(\ov{F}/F,A)$, $b \in \check{Z}^{2}(\ov{F}/F,B)$, and $x \in B(\ov{F} \otimes_{F} \ov{F})$ with $d(x) \cdot b = f(a)$, then by \cite[Construction 2.38]{Dillery1} there is a morphism $\varphi_{a,b,x} \colon \mathcal{E}_{a} \to \mathcal{E}_{b}$ of stacks over $\Spec(F)$ defined on twisted torsors by sending $(T,\psi)$ to $(T \times^{A,f}B, m_{x^{-1}} \circ \psi)$.

\subsection{The Kaletha gerbe}\label{sec:kalgerbe}

Set $u:=\varprojlim_{E/F,n} \Res_{E/F}(\mu_{n})/\mu_{n}$, where the limit is over all finite Galois extensions $E/F$ and all $n \in \mathbb{N}$, with transition maps given (for $K/E/F$ finite Galois and $m \mid n$) by the $N_{K/E}$-norm composed with the $n/m$-power map. 

We first recall \cite[\S 4.3]{Dillery1} which constructs an explicit $u$-gerbe over $\Spec(F)$ representing the class  $-1 \in H^{2}_{\tn{fppf}}(F,u) \xrightarrow{\sim} \widehat{\mathbb{Z}}$, where this isomorphism is canonical and from \cite[Theorem 3.4]{Dillery1}. Fix $(E_{k},n_{k})_{k \geq 0}$ a cofinal sequence of nested finite Galois extensions $E_{\ell}/E_{k}/F$ with $E_{0} = F$ and natural numbers $n_{k}$ such that $n_{k} \mid n_{\ell}$ if $\ell \geq k$. Denote by $\{n_{k}'\}_{k}$ the associated coprime-to-$p$ sequence and let $m_{k}$ be such that $p^{m_{k}}=n_{k}/n'_{k}$. For ease of notation, we set $u_{k}:=\Res_{E_{k}/F}(\mu_{n_{k}})/\mu_{n_{k}}$ and denote the transition map $u_{k+1} \to u_{k}$ by $\pi_{k+1,k}$. We also denote by $\pi_{k+1,k}$ the obvious extension of this transition map to a map $\Res_{E_{k+1}/F}(\mathbb{G}_{m})/\mathbb{G}_{m} \to \Res_{E_{k}/F}(\mathbb{G}_{m})/\mathbb{G}_{m}$.

We choose the special sequence of cocycles $c_{k}$ representing the canonical classes in each $H^{2}(\Gamma_{E_{k}/F},E_{k}^{\times})$ as in \cite[\S 4.4]{kaletha16}, which we identify with their corresponding \v{C}ech cocycles; as explained in \cite[\S 4.4]{kaletha16}, these cocycles are constructed via compatible sections $s_{k} \colon \Gamma_{E_{k}/F} \to W_{E_{k}/F}$  and $\zeta_{k} \colon \Gamma_{E_{k}/F} \to \Gamma_{E_{k+1}/F}$. We also choose maps $l_{k} \colon (F^{s})^{\times} \to (F^{s})^{\times}$ satisfying $l_{k}(x)^{n_{k}'} = x$ and $l_{k+1}(x)^{n'_{k+1}/n'_{k}} = l_{k}(x)$ for all $x \in (F^{s})^{\times}$---we also write $l_{k}$ for the induced map on Galois cochains. Moreover, for every $n$-cochain $y \in \mathbb{G}_{m}(\ov{F}^{\bigotimes_{F} n+1})$ with $n \geq 1$ we choose for every $k$ a $p^{m_{k}}$-root $y^{(1/p^{m_{k}})} \in \mathbb{G}_{m}(\ov{F}^{\bigotimes_{F} n+1})$ satisfying $(y^{(1/p^{m_{k+1}})})^{p^{m_{k+1}}/p^{m_{k}}} = y^{(1/p^{m_{k}})}$. We also assume that this system of roots has the following special property: If $y \in \mathbb{G}_{m}(\ov{F}^{\bigotimes_{F} n} \otimes_{F} E)$ for a finite Galois extension $E/F$, then each $y^{(1/p^{m_{k}})} \in  \mathbb{G}_{m}(\ov{F}^{\bigotimes_{F} n} \otimes_{F} E)$ as well. Such a system exists by the arguments in \cite[\S 4.3]{Dillery1}.

Denote by $\delta_{e}$ the map 
\begin{equation*}
\mu_{n_{k},F^{s}} \to u_{k,F^{s}}
\end{equation*}
induced by the map on character modules $\bbZ/n_{k}\bbZ[\Gamma_{E_{k}/F}]_{0} \to \bbZ/n_{k}\bbZ$, $f \mapsto f(e)$. Then set
\begin{equation*}
    \xi_{k} :=d[(l_{k}c_{k})^{(1/p^{m_{k}})}] \sqcup_{E_{k}/F} \delta_{e} \in \check{Z}^{2}(\ov{F}/F, u_{k})
\end{equation*}
where $``\sqcup_{E_{k}/F}$'' denotes the unbalanced cup product in \v{C}ech cohomology, defined in \cite[\S 4.2]{Dillery1}. We note that the aforementioned ``special property'' on the $p^{m_{k}}$-roots is required precisely for this unbalanced cup product to be well-defined.

Similarly, one defines
\begin{equation*}
    \alpha_{k} = [(l_{k}c_{k})^{(1/p^{m_{k}})} \sqcup_{E_{k}/F} \delta_{e,k}]^{-1} \cdot \pi_{k+1,k}( [(l_{k+1}c_{k+1})^{(1/p^{m_{k+1}})} \sqcup_{E_{k+1}/F} \delta_{e,k+1}]) \in \check{C}^{1}(\ov{F}/F,u_{k})
\end{equation*}
where $\delta_{e,k}$ denotes the obvious extension of $\delta_{e}$ to a morphism $\mathbb{G}_{m,F^{s}} \to [\Res_{E_{k}/F}(\mathbb{G}_{m})/\mathbb{G}_{m}]_{F^{s}}$.

\begin{proposition}{(\cite[Lemma 4.8]{Dillery1})}
    \begin{enumerate}
        \item{One has $d\alpha_{k} = \pi_{k+1,k}(\xi_{k+1})\xi_{k}^{-1}$.}
        \item{The element $([\xi_{k}])_{k} \in \varprojlim_{k} H_{\textnormal{fppf}}^{2}(F,u_{k}) = H_{\textnormal{fppf}}^{2}(F,u)$ equals $-1$.}
    \end{enumerate}
\end{proposition}

It follows that if we take the $u_{k}$-gerbe $\mathcal{E}_{k} := \mathcal{E}_{\xi_{k}}$ as in \cite[Definition 2.35]{Dillery1}, then one can define the $u$-gerbe
\begin{equation}\label{eq:projlimgerbe}
    \mathcal{E}:= \varprojlim_{k} \mathcal{E}_{k} \to \Spec(F),
\end{equation}
where the inverse limit is taken according to the procedure in \cite[\S 2.7]{Dillery1}, which requires the system $\{\alpha_{k}\}_{k}$. We denote the transition maps $\mathcal{E}_{k+1} \to \mathcal{E}_{k}$ by $f_{k+1,k}=\varphi_{\xi_{k+1},\xi_{k},\alpha_{k}}$. By construction, the class of this gerbe in $H_{\text{fppf}}^{2}(F,u)$ is $-1$.

\section{The main result}\label{sec:mainsec}
\subsection{The group $H^{1}(\mathcal{E}, Z \to T \to U)$}\label{sec:H1def}
This subsection is fundamentally the same as \cite[\S 5.1]{Kal26b}, but we summarize the main ideas here for completeness. We will try to emphasize only those concepts which are truly different than in loc. cit.

For the moment, fix any $u$-gerbe $\mathcal{E}$ over $F$ which is split over $\ov{F}$ and represents the class $-1 \in H_{\text{fppf}}^{2}(F,u)$ (we will specialize to the gerbe from \eqref{eq:projlimgerbe} later). 

Let $Z \to T \xrightarrow{f} U$ be a complex consisting of two tori $T$ and $U$ over $F$ and $Z$ a finite group scheme over $F$, where the first map is a monomorphism; denote by $\bar{f}$ the induced map $T/Z \to U$. We define $Z^{1}(\mathcal{E},Z \to T \to U)$ to be the set of isomorphism classes of triples $(z,c,h)$, where $z \in Z^{1}(\mathcal{E}, Z \to T)$ is as in \cite[\S 2.3]{Dillery1}, $c \in U(\ov{F})$, and $h$ is an $\ov{F}$-trivialization of the $\bar{T}:=T/Z$-torsor $\pi_{*}(z\times^{T} \bar{T})$, which sends the canonical descent datum for $\pi_{*}(z\times^{T} \bar{T})$ to a $1$-cocycle $z_{\bar{T}} \in \bar{T}(\ov{F} \otimes_{F} \ov{F})$, and we require that $\bar{f}(z_{\bar{T}}) =dc$. An isomorphism between $(z,c,h)$ and $(z',c',h')$ is an isomorphism of $T$-torsors $z \to z'$ which sends $h$ to $h'$ with $c=c'$. The set $Z^{1}(\mathcal{E}, Z \to T \to U)$ is then a group under multiplication of torsors, multiplication in $U(\ov{F})$, and sending two trivializations $h,h'$ to $h \times^{\bar{T}} h'$. For $t \in T(\ov{F})$, we abusively use $dt$ to denote the $T$-torsor with descent datum $(T_{\ov{F}}, m_{dt})$, with $\ov{F}$-trivialization $h_{\bar{t}}$ sending this descent datum to the usual \v{C}ech coboundary $d\bar{t} \in \bar{T}(\ov{F} \otimes_{F} \ov{F})$.

\begin{definition}\label{def:mainH1}
    We define $H^{1}(\mathcal{E}, Z \to T \to U)$ to be the group obtained by modding $Z^{1}(\mathcal{E},Z \to T \to U)$ out by the subgroup of ``coboundaries,'' which are all triples $(dt, f(t), h_{\bar{t}})$ for $t \in T(\ov{F})$. We define $H^{1}(\mathcal{E}_{k}, Z \to T \to U)$ for $\mathcal{E}_{k}$ as in \S \ref{sec:kalgerbe} analogously.
\end{definition}

For any other such gerbe $\mathcal{E}'$ representing $-1 \in H^{2}_{\text{fppf}}(F,u)$, there is an isomorphism of $u$-gerbes $\mathcal{E} \xrightarrow{a} \mathcal{E}'$. Moreover the proof of \cite[Lemma 2.56]{Dillery1} (using that $H^{1}_{\text{fppf}}(F,u)=0$ by \cite[Theorem 3.4]{Dillery1}) shows that, for any two such isomorphisms $a_{1}$, $a_{2}$, there is a $2$-isomorphism $\eta \colon a_{2}^{-1} \circ a_{1} \Rightarrow \text{id}_{\mathcal{E}}$ of automorphisms of $\mathcal{E}$ over $\Spec(F)$. Then, given a triple $(z,c,h)$ as above for $\mathcal{E}'$, the proof of \cite[Lemma 2.56]{Dillery1} shows that $2$-isomorphism $\eta$ (which is induced by twisting by a $1$-coboundary $dy$ for $y \in u(\ov{F})$) induces an isomorphism from $(a_{1}^{*}z,c,h)$ to $(a_{2}^{*}z,c,h)$---note that $h$ makes sense in this context because $\pi_{*}(a_{i}^{*}z \times^{T} \bar{T})$ is canonically isomorphic to $\pi'_{*}(z \times^{T} \bar{T})$, and, via these identifications, is preserved by $\eta$ because it is a $2$-morphism over $\Spec(F)$---and so one obtains the same element of $Z^{1}(\mathcal{E}, Z \to T \to U)$. It follows that the set $Z^{1}(\mathcal{E}, Z \to T \to U)$ defined above, and therefore $H^{1}(\mathcal{E}, Z \to T \to U)$ only depends on the choice of $\mathcal{E}$ up to canonical isomorphism.

\begin{remark}
    The analogous set $Z^{1}(u \to W, Z \to T \to U)$ from \cite{Kal26b} does not satisfy the above canonicity property because it is more rigid than our $Z^{1}(\mathcal{E}, Z \to T \to U)$ defined above. In cocycle language, the latter identifies pairs $(z,c)$ which differ by a $u$-coboundary twist. Since $f|_{Z}=1$, the twist by such $dz$ corresponds to multiplying a pair by a coboundary $(dz,1)$ as in the definition of coboundaries in \cite{Kal26b} (in our language, by the triple $(dz,1,\text{id}_{\bar{T}_{\ov{F}}})$), and therefore taking this modification of $Z^{1}(u \to W, Z \to T \to U)$ does not affect the ultimate cohomology group that one cares about. 
\end{remark}

When $\mathcal{E}=\mathcal{E}_{\xi}$ for $\xi \in \check{Z}^{2}(\ov{F}/F,u)$, we can use the map from \cite[Proposition 2.54]{Dillery1} to canonically identify $H^{1}(\mathcal{E}_{\xi}, Z \to T)$ with equivalence classes of $\xi$-twisted $1$-cocycles valued in $T$ with $\lambda \in \Hom_{F}(u,Z)$, which we now recall.

\begin{definition}
    For $\xi \in \check{Z}^{2}(\ov{F}/F, u)$ and a torus $T$ over $F$, a $\xi$-twisted \v{C}ech ($1$-)cocycle valued in $T$ is a pair $(y, \lambda)$, where $\lambda \in \Hom_{F}(u,T)$ and $y \in T(\ov{F} \otimes_{F} \ov{F})$ satisfies $dy = \lambda(\xi) \in T(\ov{F}^{\bigotimes_{F}3})$. We declare that $(y,\lambda)$ is cohomologous to $(y\cdot dt, \lambda)$ for all $t \in T(\ov{F})$.
\end{definition}

A fixed cocycle $v \in \check{Z}^{1}(\ov{F}/F,u)$ gives an automorphism $a_{v}:=\varphi_{\xi,\xi,v}$ of the gerbe $\mathcal{E}_{\xi}$ over $\Spec(F)$, as explained in \cite[\S 2.3]{Dillery1}. More generally, having made this identification with twisted cocycles, pullback along the map $\varphi_{a,b,v}$ described in \S \ref{sec:gerbeb} sends $(y,\lambda)$ to $(y\lambda(v),\lambda \circ f)$. In particular, pullback by $a_{dx}$ for $x \in u(\ov{F})$ on $Z^{1}(\mathcal{E}_{\xi}, Z \to T)$ sends a twisted cocycle $(y,\lambda)$ to $(\lambda(p_{2}^{\sharp}(x)p_{1}^{\sharp}(x)^{-1})y,\lambda)$.

We can extend this identification with classes of twisted cocycles from \cite[Proposition 2.54]{Dillery1} to the setting of $H^{1}(\mathcal{E}_{\xi}, Z \to T \to U)$ by sending the class of a triple $(z,c,h) \in Z^{1}(\mathcal{E}_{\xi}, Z \to T \to U)$ to the class of $(z',\lambda,c)$, where $(z',\lambda)$ is a $\xi$-twisted $1$-cocycle valued in $T$ with $\lambda \in \Hom_{F}(u,Z)$ attached to $z$ as in the proof of \cite[Proposition 2.54]{Dillery1} such that $\bar{f}(\bar{z'}) = dc$, modulo the subgroup of triples $(dt,1,f(t))$ for $t \in T(\ov{F})$. We leave the verification of the bijectivity of this map to the reader (it is very similar to the proof of Proposition 2.54 loc. cit.).

Formal arguments (along with \cite[Proposition 3.10]{Dillery1}) show that the cohomology group $H^{1}(\mathcal{E}, Z \to T \to U)$ fits into a commutative diagram completely analogous to \cite[(5.1)]{Kal26b}, replacing each rigid cohomology set with its analogue in \cite{Dillery1}, $H^{i}(\Gamma, - \to U)$ with the group $H^{i}(\ov{F}/F, - \to U)$ defined in \cite[\S A.1]{Dillery23b}, and each $H^{i}(\Gamma,Z)$ with $H^{i}_{\text{fppf}}(F,Z)$. 

Similarly, the arguments of \cite[\S A.1]{Dillery23b} show  that one has the long exact sequence
\begin{equation*}
 H^{0}(F,T) \to H^{0}(F,U) \to H^{1}(\mathcal{E}, Z \to T \to U) \to H^{1}(\mathcal{E}, Z \to T) \to H^{1}(F,U) \to
\end{equation*}
\begin{equation}\label{eq:arithLES}
    \to H^{2}(\ov{F}/F, \bar{T} \to U) \to H^{2}(F,\bar{T}) \to H^{2}(F,U) \to \dots
\end{equation}

\begin{remark}
    Since $T$ and $U$ are $F$-tori (in particular, they are smooth over $\Spec(F)$), the map
    \begin{equation*}
        H^{i}(F^{s}/F, T \to U) \to H^{i}(\ov{F}/F, T \to U)
    \end{equation*}
    is an isomorphism (see \cite[Lemma A.1]{Dillery23b}), and we may identify the left-hand term with $H^{i}(\Gamma, T(F^{s}) \to U(F^{s}))$ in the usual way.
\end{remark}

The group $H^{1}(\mathcal{E}, Z \to T)$ is always finite (it is an extension of a subgroup of $\Hom_{F}(u,Z)$ by $H^{1}(F,T)$, which are both finite), and so the image of $U(F)$ has finite index in $H^{1}(\mathcal{E}, Z \to T \to U)$. We give $H^{1}(\mathcal{E}, Z \to T \to U)$ the unique topology so that the map $U(F) \to H^{1}(\mathcal{E}, Z \to T \to U)$ from \eqref{eq:arithLES} is continuous and open.

\subsection{The group $\bar{Y}_{+,\textnormal{tor}}(Z \to T \to U)$}
Let $K/F$ be a finite Galois extension splitting $T$ and $U$, and, for a discrete $W_{K/F}$-module $M$ denote by $C_{i}(W_{K/F},M)$ the group of (inhomogeneous) finitely-supported $i$-chains $W_{K/F}^{i} \to M$, with differentials $\partial$ defined as in \cite[\S 5.1]{Kal26b}. One can then define the hyperhomology groups $H_{n}(W_{K/F}, X_{*}(T) \xrightarrow{f_{*}} X_{*}(U))$, where we use $f_{*}$ to denote the map $X_{*}(T) \to X_{*}(U)$ induced by $f$ as well as its canonical extension to a map $X_{*}(\bar{T}) \to X_{*}(U)$; recall the following explicit descriptions/related constructions from \cite[\S 5.1]{Kal26b}: 
\begin{itemize}
    \item{$Z_{0}(W_{K/F}, X_{*}(T) \to X_{*}(U)) = \{(\lambda,\mu_{1}), \lambda \in C_{0}(W_{K/F}, X_{*}(T)), \mu_{1} \in C_{1}(W_{K/F},X_{*}(U)), f_{*}(\lambda) = \partial \mu_{1}\}$.}
    \item{$B_{0}(W_{K/F}, X_{*}(T) \to X_{*}(U)) = \{(\partial \lambda_{1}, f_{*}(\lambda_{1})-\partial \mu_{2}),  \lambda_{1} \in C_{1}(W_{K/F},X_{*}(T)), \mu_{2} \in  C_{2}(W_{K/F},X_{*}(U))\}$.}
    \item{$H_{0}(W_{K/F}, X_{*}(T) \to X_{*}(U)) =Z_{0}/B_{0}$.}
    \item{$Z_{0}(W_{K/F}, X_{*}(T) \to X_{*}(U))_{0}$ is the subgroup of $Z_{0}$ of all $(\lambda,\mu_{1})$ with $N_{K/F}\lambda=0$; $H_{0}(W_{K/F}, X_{*}(T) \to X_{*}(U))_{0} = (Z_{0})_{0}/B_{0}$.}
    \item{$H_{0}(W_{K/F}, X_{*}(T) \to X_{*}(\bar{T}) \to X_{*}(U)) = Z_{0}(W_{K/F}, X_{*}(\bar{T}) \to X_{*}(U))/B_{0}(W_{K/F}, X_{*}(T) \to X_{*}(U)).$}  
     \item{$H_{0}(W_{K/F}, X_{*}(T) \to X_{*}(\bar{T}) \to X_{*}(U))_{0} = Z_{0}(W_{K/F}, X_{*}(\bar{T}) \to X_{*}(U))_{0}/B_{0}(W_{K/F}, X_{*}(T) \to X_{*}(U)).$}
    
\end{itemize}
For $L/K/F$ finite Galois, there is a coinflation map \cite[\S 5.1]{Kal26b} $C_{n}(W_{L/F},X_{*}(T)) \to C_{n}(W_{K/F},X_{*}(T))$ inducing maps on hyperhomology $H_{n}(W_{L/F}, X_{*}(T) \to X_{*}(U)) \to H_{n}(W_{K/F}, X_{*}(T) \to X_{*}(U))$. 

All of the above groups satisfy the same properties and compatibilities, namely Facts 5.1-5.5 in \cite{Kal26b}, as their mixed-characteristic analogues. The arguments are unchanged because all calculations only use \cite[\S A.3]{KS99}, which is valid over a local function field (cf. \cite[\S A.2]{Dillery23b}). In particular, Fact 5.4 holds, which says that coinflation induces an isomorphism 
\begin{equation*}
    H_{0}(W_{L/F}, X_{*}(T) \to X_{*}(\bar{T}) \to X_{*}(U))_{0} \to H_{0}(W_{K/F}, X_{*}(T) \to X_{*}(\bar{T}) \to X_{*}(U))_{0}
\end{equation*}
for all $L/K/F$.

Recall that $\bar{Y}_{+,\text{tor}}(Z \to T)$ is defined as $\varinjlim_{k} X_{*}(\bar{T})^{N_{E_{k}/F}}/I_{\Gamma_{E_{k}/F}}X_{*}(T)$, where the superscript ``$N_{E_{k}/F}$'' means all of those elements killed by the $\Gamma_{E_{k}/F}$-norm map, $I_{\Gamma_{E_{k}/F}}$ denotes the augmentation ideal in $\mathbb{Z}[\Gamma_{E_{k}/F}]$, and $Y_{\text{tor}}(U):=X_{*}(U)_{\Gamma,\text{tor}}$.

\begin{definition}
    We define $\bar{Y}_{+,\textnormal{tor}}(Z \to T \to U)$ as the inverse limit over all finite Galois $L/F$ (splitting $T$ and $U$) of the groups $H_{0}(W_{L/F}, X_{*}(T) \to X_{*}(\bar{T}) \to X_{*}(U))_{0}$, where the transition maps are induced by coinflation.
\end{definition}

As in \cite[Fact 5.5]{Kal26b}, there is an exact sequence
\begin{equation}\label{eq:dualLES}
   H_{1}(W_{F},X_{*}(T)):= \varprojlim_{L, \text{coinf}}H_{1}(W_{L/F}, X_{*}(T)) \to H_{1}(W_{F}, X_{*}(U)) \to 
   \end{equation}
   \begin{equation*}
   \to \bar{Y}_{+,\textnormal{tor}}(Z \to T \to U) \to \bar{Y}_{+,\textnormal{tor}}(Z \to T) \to Y_{\textnormal{tor}}(U).
\end{equation*}
Langlands duality \cite{Borel77} (valid for non-archimedean local fields of any characteristic) says that the above projective limits stabilize and are isomorphic to $T(F)$ and $U(F)$, respectively, and we give them the induced topology. We give $\bar{Y}_{+,\textnormal{tor}}(Z \to T \to U)$ the unique topology so that the map from $H_{1}(W_{F}, X_{*}(U))$ is continuous and open.

The dual morphism $\hat{f} \colon \widehat{U} \to \widehat{T}$ lifts uniquely to $\hat{\bar{f}} \colon \widehat{U} \to \widehat{\bar{T}}$. Setting $\widehat{Z}:=\text{Ker}[\widehat{\bar{T}} \to \widehat{T}]$ and $[\widehat{\bar{T}}]^{\Gamma,+}$ the preimage of $\widehat{T}^{\Gamma}$ in $\widehat{\bar{T}}$, we define $Z^{1}_{\text{cts}}(W_{F}, \widehat{Z} \to \widehat{\bar{T}} \leftarrow \widehat{U})$ to be the group of pairs $(z,\dot{c})$ where $z \in Z^{1}_{\cts}(W_{F},\widehat{U})$ and $\dot{c} \in \widehat{\bar{T}}$ such that $\hat{f}(z)=dc$, where $c$ is the image of $\dot{c}$ in $\widehat{T}$, $B^{1}(W_{F}, \widehat{Z} \to \widehat{\bar{T}} \leftarrow \widehat{U})$ the subgroup of all pairs $(du, \hat{\bar{f}}(u))$ for $u \in \widehat{U}$, and $H^{1}_{\text{cts}}(W_{F}, \widehat{Z} \to \widehat{\bar{T}} \leftarrow \widehat{U}) :=Z^{1}_{\text{cts}}(W_{F}, \widehat{Z} \to \widehat{\bar{T}} \leftarrow \widehat{U})/B^{1}(W_{F}, \widehat{Z} \to \widehat{\bar{T}} \leftarrow \widehat{U})$. One has an exact sequence
\begin{equation*}
    H^{1}_{\cts}(W_{F}, \widehat{T}) \leftarrow  H^{1}_{\cts}(W_{F}, \widehat{U}) \leftarrow H^{1}_{\text{cts}}(W_{F}, \widehat{Z} \to \widehat{\bar{T}} \leftarrow \widehat{U}) \leftarrow [\widehat{\bar{T}}]^{\Gamma,+} \leftarrow \widehat{U}^{\Gamma}
\end{equation*}
and we define $H^{1}_{\text{cts}}(W_{F}, \widehat{Z} \to \widehat{\bar{T}} \leftarrow \widehat{U})_{\text{red}}$ as the quotient of $H^{1}_{\text{cts}}(W_{F}, \widehat{Z} \to \widehat{\bar{T}} \leftarrow \widehat{U})$ by the image of $[\widehat{\bar{T}}]^{\Gamma,+,\circ}$.

\subsection{Duality on the Galois side}
As shown at the beginning of \cite[\S 5.2]{Kal26b}, one has a pairing sending $(z,\dot{c}) \in Z^{1}_{\text{cts}}(W_{F}, \widehat{Z} \to \widehat{\bar{T}} \leftarrow \widehat{U})$ and $(\bar{\lambda},\mu_{1}) \in  Z_{0}(W_{K/F}, X_{*}(\bar{T}) \to X_{*}(U))_{0}$ to the value
\begin{equation*}
    \langle (z,\dot{c}) , (\bar{\lambda},\mu_{1}) \rangle_{K} := \langle \dot{c}, \bar{\lambda} \rangle_{\bar{T}} \cdot \prod_{w \in W_{K/F}} \langle z(w), \mu_{1}(w) \rangle_{U}^{-1} \in \mathbb{C}^{\times},
\end{equation*}
where the pairings used in the above formula are the usual character-pairings for $\bar{T}$ and $U$, respectively. Moreover, the above pairing is compatible with coinflation, kills the appropriate coboundaries (and the image of $[\widehat{\bar{T}}]^{\Gamma,+,\circ}$), and induces a pairing
\begin{equation}\label{eq:dualpairing}
H^{1}_{\text{cts}}(W_{F}, \widehat{Z} \to \widehat{\bar{T}} \leftarrow \widehat{U})_{\text{red}} \times \bar{Y}_{+,\textnormal{tor}}(Z \to T \to U) \to \mathbb{C}^{\times}
\end{equation}
which is functorial in $Z \to T \to U$.

\begin{proposition}\label{prop:pp1}
    The above pairing is compatible with the pairing $ \pi_{0}([\widehat{\bar{T}}]^{\Gamma,+}) \times \bar{Y}_{+,\textnormal{tor}}(Z \to T) \to \mathbb{C}^{\times}$ from \cite[Proposition 7.10]{Dillery1}, with the negative of Langlands' pairing $H_{\cts}^{1}(W_{F},\widehat{U}) \times H_{1}(W_{F},X_{*}(U))\to \mathbb{C}^{\times}$, and induces an isomorphism
    \begin{equation*}
       H^{1}_{\textnormal{cts}}(W_{F}, \widehat{Z} \to \widehat{\bar{T}} \leftarrow \widehat{U})_{\textnormal{red}} \to \Hom_{\cts}(\bar{Y}_{+,\textnormal{tor}}(Z \to T \to U),\mathbb{C}^{\times}).
    \end{equation*}
\end{proposition}

\begin{proof}
    The proof of the analogous mixed characteristic statement \cite[Fact 5.6]{Kal26b} goes through verbatim here, replacing the use of \cite[Proposition 5.3]{kaletha16} with \cite[Proposition 7.10]{Dillery1}.
\end{proof}

\subsection{Arithmetic duality}
We use the objects and notation introduced in \S \ref{sec:setup}, and continue with a fixed complex $Z \to T \xrightarrow{f} U$ as in the previous subsections. Fix $k \gg 0$ such that $T$ and $U$ are split over $E_{k}$ and $|Z|$ divides $n_{k}$. We now fix the $u$-gerbe $\mathcal{E}$ defined in \eqref{eq:projlimgerbe} representing $-1 \in H^{2}(F,u)$.

Fix $(\bar{\lambda},\mu_{1}) \in Z_{0}(W_{E_{k}/F}, X_{*}(\bar{T}) \to X_{*}(U))_{0}$. Recall from \cite[\S A.3]{KS99} the element
\begin{equation*}
    \phi_{U}(\mu_{1}):=\prod_{(\sigma,\tau,a) \in \Gamma_{E_{k}/F} \times \Gamma_{E_{k}/F} \times E_{k}^{\times}} \sigma(\mu_{1}(as_{k}(\tau)))(c_{k}(\sigma,\tau)^{-1}\sigma(a)^{-1}) \in U(E_{k}).
\end{equation*}
As in \cite[\S 5.2]{Kal26b}, we define $\dot{\psi}_{\bar{T}}(\bar{\lambda}) \in Z^{1}(\mathcal{E},Z \to T)$ to be the pullback to $\mathcal{E}$ along the canonical projection $\mathcal{E} \to \mathcal{E}_{k}$ of $Z_{\bar{\lambda},k} \in Z^{1}(\mathcal{E}_{k}, Z \to T)$ as defined in \cite[\S 4.3, before Lemma 4.9]{Dillery1}, which corresponds to the $\xi_{k}$-twisted $1$-cocycle $(z_{\bar{\lambda},k}, \phi_{\bar{\lambda},k})$ where $z_{\bar{\lambda},k}= (l_{k}c_{k})^{(1/p^{m_{k}})} \sqcup_{E_{k}/F} n_{k}\bar{\lambda} \in T(\ov{F} \otimes_{F} \ov{F})$ and $\phi_{\bar{\lambda},k}$ is the image of $\bar{\lambda}$ along the composition $[X_{*}(\bar{T})/X_{*}(T)]^{N_{E_{k}/F}} \to \Hom_{E_{k}}(\mu_{n_{k}},Z)^{N_{E_{k}/F}} \xrightarrow{\sim} \Hom_{F}(u_{k},Z)$.

The $\bar{T}$-torsor (which can be computed at the $\mathcal{E}_{k}$-level)
\begin{equation*}
\pi_{*}(Z_{\bar{\lambda},k} \times^{T} \bar{T})
\end{equation*}
is given by the descent datum $(\bar{T}_{\ov{F}}, \bar{z}_{\bar{\lambda},k})$ and so we can take the structure trivialization $h$ of $\pi_{*}(Z_{\bar{\lambda},k} \times^{T} \bar{T})$ and obtain the 1-cocycle $\bar{z}_{\bar{\lambda},k} \in \bar{T}(\ov{F} \otimes_{F} \ov{F})$. Since $\bar{\lambda} \in X_{*}(\bar{T})$ we can move $n_{k}$ to the other side of the cup product and rewrite $\bar{z}_{\bar{\lambda},k}$ as $c_{k} \cup \bar{\lambda} = \psi_{\bar{T}}(\bar{\lambda})$ (where the map $\psi_{\bar{T}}$ is as defined in \cite[\S A.3]{KS99}). The functoriality of $\psi_{-}$ then implies (see \cite[\S 5.2]{Kal26b}) that the triple $(\dot{\psi}_{\bar{T}}(\bar{\lambda}), \phi_{U}(\mu_{1}),h)$ gives an element of $Z^{1}(\mathcal{E},Z \to T \to U)$.

Moreover, if one assumes that $(\bar{\lambda},\mu_{1}) = (\partial\lambda_{1}, f_{*}(\lambda_{1}) - \partial\mu_{2}) \in B_{0}(W_{E_{k}/F}, X_{*}(T) \to X_{*}(U))$ then $\partial \lambda_{1} \in X_{*}(T)$ implies that $\dot{\psi}_{\bar{T}}(\partial \lambda_{1}) = \psi_{T}(\partial \lambda_{1}) = d \phi_{T}(\lambda_{1})$ (using \cite[\S A.3]{KS99} for this last equality) and $\phi_{U}(f_{*}(\lambda_{1})-\partial \mu_{2}) = f(\phi_{T}(\lambda_{1}))$, and so $(\partial\lambda_{1}, f_{*}(\lambda_{1})-\partial\mu_{2})$ maps via $\dot{\psi}_{\bar{T}} \times \phi_{U}$ to a coboundary as in Definition \ref{def:mainH1}. We therefore have a well-defined group homomorphism
\begin{equation*}
    H_{0}(W_{E_{k}/F}, X_{*}(T) \to X_{*}(\bar{T}) \to X_{*}(U))_{0} \xrightarrow{\dot{\psi}_{\bar{T}} \times \phi_{U}} H^{1}(\mathcal{E},Z \to T \to U).
\end{equation*}

For $(\bar{\lambda},\mu_{1})$ as above at level $k+1$, \cite[Lemma 4.9]{Dillery1} says that the pullback of $Z_{\bar{\lambda},k}$ along $\mathcal{E}_{k+1} \to \mathcal{E}_{k}$ is $Z_{\bar{\lambda},k+1}$, and hence $\dot{\psi}_{\bar{T},k}(\text{coinf}(\bar{\lambda})) = \dot{\psi}_{\bar{T},k+1}(\bar{\lambda})$. With this equality in hand, the same argument as in \cite[\S 5.2, before (5.5)]{Kal26b} implies that the maps $\dot{\psi}_{\bar{T},k} \times \phi_{U,k}$ are compatible across the coinflation maps for varying $k$ (even at the level of cocycles), and therefore give a homomorphism
\begin{equation}\label{eq:TNmap}
    \bar{Y}_{+,\textnormal{tor}}(Z \to T \to U) \to H^{1}(\mathcal{E},Z \to T \to U).
\end{equation}

The following is a direct analogue of \cite[Proposition 5.7]{Kal26b} and is our main result:
\begin{theorem}\label{thm:TN}
    The map \eqref{eq:TNmap} is a functorial isomorphism which does not depend on the choices made during its construction.
\end{theorem}

\begin{proof}
The map \eqref{eq:TNmap} is functorial because the maps $\dot{\psi}_{\bar{T},k} \times \phi_{U,k}$ are functorial. The fact that it is an isomorphism comes from applying the five lemma to the exact sequences \eqref{eq:arithLES} and \eqref{eq:dualLES}, using Langlands duality for tori, Tate--Nakayama duality, and \cite[Theorem 4.10]{Dillery1} to supply isomorphisms for all other terms in the sequences. 

As is the case in \cite{Kal26b}, more work is required to show that \eqref{eq:TNmap} does not depend on the choice of cocycles $c_{k}$---which themselves are determined by the sections $s_{k}$, $\zeta_{k}$---and the root maps $l_{k}$, $(-)^{(1/p^{m_{k}})}$.

Take any other $s_{k}'$, $\zeta_{k}'$, $l_{k}'$, and $(-)^{(1/p^{m_{k}})'}$ satisfying the properties of \cite[\S 4.4]{kaletha16} and \S \ref{sec:setup}, giving a different gerbe $\mathcal{E}' = \varprojlim_{k} \mathcal{E}'_{k}$ and map $ \bar{Y}_{+,\textnormal{tor}}(Z \to T \to U) \to H^{1}(\mathcal{E}', Z \to T \to U)$. Recall that $\mathcal{E}$ is isomorphic to $\mathcal{E}'$, and, by the discussion in \S \ref{sec:H1def}, all choices of isomorphism induce the same map $H^{1}(\mathcal{E}, Z \to T \to U) \to  H^{1}(\mathcal{E}', Z \to T \to U)$. Recall from \cite{Kal26b} the map $\eta_{k} \colon \Gamma_{E_{k}/F} \to E_{k}^{\times}$ defined by the relation $s'_{k}(\sigma)=\eta_{k}(\sigma)s_{k}(\sigma)$ and the $u_{k}$-valued $1$-cochain
\begin{equation*}
\alpha_{k',k}:=[(l_{k}'c_{k}')^{(1/p^{m_{k}})'} \cdot ((l_{k}c_{k})^{(1/p^{m_{k}})})^{-1} \cdot (d((l_{k}\eta_{k})^{(1/p^{m_{k}}))})^{-1}] \sqcup_{E_{k}/F} \delta_{e}.
\end{equation*}

This cochain defines an isomorphism of gerbes $\bar{g}_{k} \colon \mathcal{E}'_{k} \to \mathcal{E}_{k}$ as follows. It sends a $\xi_{k}'$-twisted $u_{k}$-torsor $(T,\psi)$ over an $F$-scheme $S$ to $(T,m_{\alpha_{k',k}^{-1}} \circ \psi)$, where $m_{\alpha_{k',k}}$ denotes the automorphism of the $u_{k,S_{\ov{F} \otimes_{F} \ov{F}}}$-torsor $T_{S_{\ov{F} \otimes_{F} \ov{F}}}$ given by acting by the image of $\alpha_{k',k}$ in $u_{k}(S_{\ov{F} \otimes_{F} \ov{F}})$. We claim that the above formula gives a well-defined map, and that pullback by this map sends each $T$-torsor $Z_{\bar{\lambda},k}$ (where $\bar{\lambda} \in X_{*}(\bar{T})^{N_{E_{k}/F}}$) to the $T$-torsor on $\mathcal{E}_{k}'$ corresponding to the $\xi_{k}'$-twisted $1$-cocycle $(z_{\bar{\lambda},k}' \cdot d((l_{k}\eta_{k})^{(1/p^{m_{k}})} \sqcup_{E_{k}/F} n_{k}\bar{\lambda})^{-1}, \phi_{\bar{\lambda},k})$. To see the first part, we compute (using \cite[Proposition 4.4]{Dillery1}) that
\begin{equation*}
    d\alpha_{k',k} = \xi_{k}'\xi_{k}^{-1}.
\end{equation*}
The second part is immediate from the $\xi_{k}$-twisted $1$-cocycle description of $Z_{\bar{\lambda},k}$.

These comparison maps $\mathcal{E}'_{k} \xrightarrow{\bar{g}_{k}} \mathcal{E}_{k}$ are not naively compatible for varying $k$, so one introduces the correction factor $1$-cochain $\beta_{k}:=[(l_{k}\eta_{k})^{(1/p^{m_{k}})}]^{-1}N_{E_{k+1}/E_{k}}((l_{k}\eta_{k+1})^{(1/p^{m_{k}})}) \in \mathbb{G}_{m}(\ov{F} \otimes_{F} E_{k}) = \Res_{E_{k}/F}(\mathbb{G}_{m})(\ov{F})$, where $N_{E_{k+1}/E_{k}}$ is the map $\mathbb{G}_{m}(\ov{F} \otimes_{F} E_{k+1}) \to \mathbb{G}_{m}(\ov{F} \otimes_{F} E_{k})$ given by the composition
\begin{equation*}
  \mathbb{G}_{m}(\ov{F} \otimes_{F} E_{k+1}) \xrightarrow{\sim}  \prod_{\Gamma_{E_{k+1}/F}} \ov{F}^{\times} \to \prod_{\Gamma_{E_{k}/F}} \ov{F}^{\times} \xrightarrow{\sim} \mathbb{G}_{m}(\ov{F} \otimes_{F} E_{k}),
\end{equation*}
where the middle map sends a function $f$ to $[a \mapsto \prod_{b \mapsto a} f(b)]$.

The key result concerning $\beta_{k}$ is our analogue of \cite[Lemma 5.9]{Kal26b}, which says that the image of $\beta_{k}^{n_{k}}$ via the map $\Res_{E_{k}/F}(\mathbb{G}_{m})(\ov{F}) \to [\Res_{E_{k}/F}(\mathbb{G}_{m})/\mathbb{G}_{m}](\ov{F})$ is trivial (implying that $\beta_{k} \in u_{k}(\ov{F})$) and $f_{k+1,k} \circ \bar{g}_{k+1} = a_{d\beta_{k}^{-1}} \circ \bar{g}_{k} \circ f_{k+1,k}'$. Recall from \S \ref{sec:H1def} that $a_{d\beta_{k}^{-1}}$ denotes the automorphism of $\mathcal{E}_{k}$ induced by the $1$-coboundary of $\beta_{k}^{-1} \in u_{k}(\ov{F})$. As in \cite{Kal26b}, we directly compute on $\xi_{k+1}'$-twisted $u_{k+1}$-torsors $(T,\psi)$ that the image of $(T,\psi)$ under
\begin{equation*}
    f_{k+1,k} \circ \bar{g}_{k+1}
\end{equation*}
differs from its image under $\bar{g}_{k} \circ f'_{k+1,k}$ by multiplying the resulting $\xi_{k}$-twisted descent datum isomorphism by the inverse of
\begin{equation*}
\pi_{k+1,k}(d((l_{k+1}\eta_{k+1})^{(1/p^{m_{k+1}})}) \sqcup_{E_{k+1}/F} \delta_{e,k+1})^{-1} \cdot (d((l_{k}\eta_{k})^{(1/p^{m_{k}})}) \sqcup_{E_{k}/F} \delta_{e,k}) \in u_{k}(\ov{F} \otimes_{F} \ov{F}),
\end{equation*}
(in the language of \S \ref{sec:gerbeb}, the resulting gerbe automorphism is $\varphi_{\xi_{k},\xi_{k},x}$, where $x$ is the above cochain) and we can pull the differential outside (again using \cite[Proposition 4.4]{Dillery1} and $d\delta_{e,k}=0$) to rewrite the above expression as
\begin{equation}\label{eq:dstg}
d[\pi_{k+1,k}((l_{k+1}\eta_{k+1})^{(1/p^{m_{k+1}})} \sqcup_{E_{k+1}/F} \delta_{e,k+1})^{-1} \cdot ((l_{k}\eta_{k})^{(1/p^{m_{k}})} \sqcup_{E_{k}/F} \delta_{e,k})].
\end{equation}

The element $(l_{k}\eta_{k})^{(1/p^{m_{k}})} \sqcup_{E_{k}/F} \delta_{e,k}$ lies in $[\Res_{E_{k}/F}(\mathbb{G}_{m})/\mathbb{G}_{m}](\ov{F})$; one computes explicitly (using the unbalanced cup product formula in \cite[\S 4.2]{Dillery1} and Lemma 4.5 loc. cit.) that 
\begin{equation*}
    (l_{k}\eta_{k})^{(1/p^{m_{k}})} \sqcup_{E_{k}/F} \delta_{e,k} = \prod_{a \in \Gamma_{E_{k}/F}} ({^{a}\delta}_{e,k})[(l_{k}\eta_{k})^{(1/p^{m_{k}})}(a)]
\end{equation*}
(where for $x \in \mathbb{G}_{m}(\ov{F} \otimes_{F} E_{k})$ we denote its image in $\prod_{a \in \Gamma_{E_{k}/F}} \ov{F}^{\times}$ by $(x(a))_{a}$), which is exactly the image of $(l_{k}\eta_{k})^{(1/p^{m_{k}})} \in \mathbb{G}_{m}(\ov{F} \otimes_{F} E_{k}) = \mathrm{Res}_{E_{k/F}}(\mathbb{G}_{m})(\ov{F})$ in $[\Res_{E_{k}/F}(\mathbb{G}_{m})/\mathbb{G}_{m}](\ov{F})$. Applying this formula for $k+1$ as well reduces \eqref{eq:dstg} to $d\beta_{k}^{-1}$, as claimed. The proof that $\beta_{k}^{n_{k}}=1$ in $[\Res_{E_{k}/F}(\mathbb{G}_{m})/\mathbb{G}_{m}](\ov{F})$ is then a Galois cocycle computation identical to the one in the proof of \cite[Lemma 5.9]{Kal26b}.

Putting everything together, we choose $\dot{\beta}_{k} \in u(\ov{F})$ mapping to $\beta_{k}$ under $u \xrightarrow{\pi_{k}} u_{k}$, set $\dot{\beta}_{<k} := \prod_{i=0}^{k-1} \dot{\beta}_{i}$, and define $g_{k} \colon \mathcal{E}_{k}' \to \mathcal{E}_{k}$ as $a_{d\pi_{k}(\dot{\beta}_{<k})} \circ \bar{g}_{k}$. By construction, the system of isomorphisms $(g_{k})_{k}$ is compatible with the transition maps $f_{k+1,k}$ and $f'_{k+1,k}$ for all $k$ and therefore induces an isomorphism of gerbes $g \colon \mathcal{E}' \to \mathcal{E}$. Our above work shows that the pullback of any $Z_{\bar{\lambda},k}$ by $g_{k}$, denoted by $Z_{\bar{\lambda},k}''$, corresponds to the $\xi_{k}'$-twisted cocycle
\begin{equation*}
    (z''_{\bar{\lambda},k}, \phi_{\bar{\lambda},k}):=(z'_{\bar{\lambda},k} \cdot \phi_{\bar{\lambda},k}(\pi_{k}(p_{2}^{\sharp}(\dot{\beta}_{<k})p_{1}^{\sharp}(\dot{\beta}_{<k})^{-1})) \cdot d((l_{k}\eta_{k})^{(1/p^{m_{k}})} \sqcup_{E_{k}/F} n_{k}\bar{\lambda})^{-1},\phi_{\bar{\lambda},k}),
\end{equation*}
which equals (because $\phi_{\bar{\lambda},k}$ is defined over $F$)
\begin{equation*}
     (z'_{\bar{\lambda},k} \cdot d(\phi_{\bar{\lambda},k}(\pi_{k}(\dot{\beta}_{<k}))^{-1} \cdot ((l_{k}\eta_{k})^{(1/p^{m_{k}})} \sqcup_{E_{k}/F} n_{k}\bar{\lambda}))^{-1},\phi_{\bar{\lambda},k}).
\end{equation*}

We conclude that for $(\bar{\lambda},\mu_{1}) \in Z_{0}(W_{E_{k}/F},X_{*}(\bar{T}) \to X_{*}(U))_{0}$, mapping to $Z^{1}(\mathcal{E}_{k}, Z \to T \to U)$ via $\dot{\psi}_{\bar{T},k} \times \phi_{U,k}$ and then pulling back by $g_{k}$ yields 
\begin{equation*}
(Z''_{\bar{\lambda},k}, \phi_{U,k}(\mu_{1}),h)
\end{equation*}
which differs from $(Z_{\bar{\lambda},k}', \phi'_{U,k}(\mu_{1}),h)$ by the coboundary (written in terms of twisted cocycles using \S \ref{sec:H1def})
\begin{equation*}
(d(\phi_{\bar{\lambda},k}(\pi_{k}(\dot{\beta}_{<k}))^{-1} \cdot ((l_{k}\eta_{k})^{(1/p^{m_{k}})} \sqcup_{E_{k}/F} n_{k}\bar{\lambda}))^{-1},1,[\eta_{k}\cup f_{*}(\bar{\lambda})]^{-1}),
\end{equation*}
using the identity $\phi_{U,k}(\mu_{1}) = \phi'_{U,k}(\mu_{1}) \cdot(\eta_{k} \cup f_{*}(\bar{\lambda}))^{-1}$ from \cite[\S A.3]{KS99}. We have also used the fact that $f|_{Z}=1$ and so the value $f(\phi_{\bar{\lambda},k}(\pi_{k}(\dot{\beta}_{<k}))^{-1} \cdot ((l_{k}\eta_{k})^{(1/p^{m_{k}})} \sqcup_{E_{k}/F} n_{k}\bar{\lambda}))$ may be computed after applying the map $T \to \bar{T}$ which kills $\phi_{\bar{\lambda},k}(\pi_{k}(\dot{\beta}_{<k}))$ and allows one to pull the $n_{k}$ across the cup product, reducing the expression to $\bar{f}(\eta_{k} \cup \bar{\lambda}) = \eta_{k} \cup f_{*}(\bar{\lambda})$. This finishes the proof that the map \eqref{eq:TNmap} is independent of $s_{k}$, $\zeta_{k}$, $l_{k}$, and $(-)^{(1/p^{m_{k}})}$. From here, the proof of independence of \eqref{eq:TNmap} on the choice of $(E_{k},n_{k})$ is identical to that of \cite{Kal26b}.
\end{proof}

Combining the isomorphism \eqref{eq:TNmap} with the pairing \eqref{eq:dualpairing} gives a pairing
\begin{equation}\label{eq:TNbis}
    H^{1}(\mathcal{E}, Z \to T\to U) \times H^{1}_{\textnormal{cts}}(W_{F}, \widehat{Z} \to \widehat{\bar{T}} \leftarrow \widehat{U})_{\textnormal{red}} \to \mathbb{C}^{\times}
\end{equation}
which, using Proposition \ref{prop:pp1} and Corollary \ref{cor:TN} below, is a perfect pairing.

The isomorphism \eqref{eq:TNmap} and pairing \eqref{eq:TNbis} satisfy the following compatibilities, whose proofs, given what we have already established above, are unchanged from their mixed-characteristic analogues \cite[Lemma 5.10, Corollary 5.11]{Kal26b} which we therefore omit:

\begin{corollary}\label{cor:TN}
    \begin{enumerate}
        \item The composition $H_{1}(W_{F},X_{*}(U)) \to  \bar{Y}_{+,\textnormal{tor}}(Z \to T \to U) \xrightarrow{\eqref{eq:TNmap}} H^{1}(\mathcal{E}, Z \to T \to U)$ is the negative of the composition $H_{1}(W_{F}, X_{*}(U)) \xrightarrow{\sim} H^{0}(F,U) \to H^{1}(\mathcal{E}, Z \to T \to U)$, where the first isomorphism is Langlands' duality for tori.
        \item The composition $\bar{Y}_{+,\textnormal{tor}}(Z \to T \to U) \xrightarrow{\eqref{eq:TNmap}} H^{1}(\mathcal{E}, Z \to T \to U) \to H^{1}(\mathcal{E}, Z \to T)$ agrees with the composition $\bar{Y}_{+,\textnormal{tor}}(Z \to T \to U) \to \bar{Y}_{+,\textnormal{tor}}(Z \to T) \xrightarrow{\sim} H^{1}(\mathcal{E}, Z \to T)$, where the last isomorphism is \cite[Theorem 4.10]{Dillery1}.
        \item The pairing \eqref{eq:TNbis} is compatible with Langlands' pairing $H^{0}(F,U) \times H^{1}_{\cts}(W_{F}, \widehat{U}) \to \mathbb{C}^{\times}$. 
        \item The pairing \eqref{eq:TNbis} is compatible with the pairing $H^{1}(\mathcal{E}, Z \to T) \times [\widehat{\bar{T}}]^{\Gamma,+} \to \mathbb{C}^{\times}$ from \cite[Corollary 7.11]{Dillery1}.
    \end{enumerate}
\end{corollary}

\section{Applications}
We briefly recall the setup of \cite[\S 8.3]{Kal26b}. Fix an $F$-rational torus $T$, $A$ a finite \'{e}tale group scheme over $F$ with action on $T$ defined over $F$, $Z \subset T$ a finite $F$-rational subgroup on which $A$ acts trivially, $\phi \colon W_{F} \to \widehat{T} \rtimes W_{F}$ and $z \in Z^{1}(\mathcal{E}, Z \to T)$ along with an $\ov{F}$-trivialization $h$ of $\bar{z}:=\pi_{*}(z \times^{T} \bar{T})$. Suppose we have two extensions (the left-hand extension being of group schemes over $F$, the right-hand being of abstract groups with a $\Gamma$-action)
\begin{equation*}
    1 \to T \to \tilde{T} \to A \to 1, \hspace{1mm} 1 \to \widehat{\bar{T}} \to \bar{\mathcal{T}} \to A(F^{s}) \to 1;
\end{equation*}
the right-hand extension induces, via modding out by $\widehat{Z}$, the extension $1 \to \widehat{T} \to \mathcal{T} \to A(F^{s}) \to 1$. Assume that the following two induced sequences are also exact (setting $A^{\Gamma}:= A(F^{s})^{\Gamma} = A(F)$):
\begin{equation*}
    1 \to T(F) \to \tilde{T}(F) \to A^{\Gamma} \to 1, \hspace{1mm} 1 \to \widehat{T}^{\Gamma} \to \mathcal{T}^{\Gamma} \to A^{\Gamma} \to 1.
\end{equation*}

Since $A$ acts $F$-rationally on $T$, $A^{\Gamma}$ acts on the set of $\widehat{T}$-conjugacy classes of homomorphisms $\phi \colon W_{F} \to \widehat{T} \rtimes W_{F}$ and also on $H^{1}(\mathcal{E}, Z \to T)$; we denote by $A^{[\phi]}$ and $A^{[z]}$ the stabilizer in $A^{\Gamma}$ of $[\phi]$ or $[z]$ (respectively). Our assumptions imply that we have exact sequences 
\begin{equation}\label{eq:invSES}
1 \to T(F) \to \tilde{T}(F)^{[\phi]} \to A^{[\phi]} \to 1, \hspace{1mm} 1 \to [\widehat{\bar{T}}]^{\Gamma,+} \to \bar{\mathcal{T}}^{\Gamma,+,[z]} \to A^{[z]} \to 1.
\end{equation}
In the above expression, we denote by $[\widehat{\bar{T}}]^{\Gamma,+}$, resp. $\bar{\mathcal{T}}^{\Gamma,+}$, the preimage of $\widehat{T}^{\Gamma}$, resp. $\mathcal{T}^{\Gamma}$, in $\widehat{\bar{T}}$, resp. $\bar{\mathcal{T}}$, by $\tilde{T}(F)^{[\phi]}$ the preimage of $A^{[\phi]}$ in $\tilde{T}(F)$, and by $\bar{\mathcal{T}}^{\Gamma,+,[z]}$ the preimage of $A^{[z]}$ in $\bar{\mathcal{T}}^{\Gamma,+}$. Langlands duality lets us view $[\phi]$ as a character $T(F) \xrightarrow{[\phi]} \mathbb{C}^{\times}$, and \cite[Corollary 7.11]{Dillery1} lets us view $[z]$ as a character $[\widehat{\bar{T}}]^{\Gamma,+} \xrightarrow{[z]} \mathbb{C}^{\times}$. We denote by $\mathcal{E}_{[\phi]}^{0}$ (resp. $\mathcal{E}_{[z]}^{0}$) the push-out of the left-hand (resp. right-hand) sequence of \eqref{eq:invSES} by $[\phi]$ (resp. by $[z]$).

As in \cite[\S 3.2]{Kal26b}, one defines the inner form $\tilde{T}_{z}$ with $\tilde{T}_{z}(F) = \{\tilde{t} \in \tilde{T}(\ov{F}) \mid p_{1}^{\sharp}(\tilde{t})=\text{Ad}(\bar{z})p_{2}^{\sharp}(\tilde{t}) \}$, where $\bar{z}$ abusively denotes the image of $z$ in $\check{Z}^{1}(\ov{F}/F,\bar{T})$ obtained from the $\bar{T}$-torsor $\bar{z}$ and $h$. We obtain an extension
\begin{equation*}
1 \to T(F) \to \tilde{T}_{z}(F) \to A^{[z]} \to 1. \end{equation*}
One then pulls this extension back along $A^{[z]} \cap A^{[\phi]}=:A^{[z],[\phi]} \to A^{[z]}$ and then takes the push-out along $[\phi]$ to obtain $\mathcal{E}^{z}_{[\phi]}$, an extension of $A^{[z],[\phi]}$ by $\mathbb{C}^{\times}$.

The dual counterpart of the above extension is obtained by setting $S_{\phi} := \mathcal{T}^{\phi(W_{F})}$, which surjects onto $A^{[\phi]}$, with preimage $S_{\phi}^{+}$ in $\bar{\mathcal{T}}$, which is an extension of $A^{[\phi]}$ by $[\widehat{\bar{T}}]^{\Gamma,+}$. One pulls this back along $A^{[z],[\phi]} \to A^{[\phi]}$ and takes the push-out by $[z]$ to obtain the extension $\mathcal{E}_{[z]}^{\phi}$.

The following is the function field analogue of \cite[Proposition 8.2]{Kal26b}:

\begin{proposition}\label{prop:prop8.2}
  Denote by $\mathcal{E}_{[\phi]}^{0,[z]}$ and $\mathcal{E}_{[z]}^{0,[\phi]}$ the pullbacks of $\mathcal{E}_{[\phi]}^{0}$ and $\mathcal{E}_{[z]}^{0}$ along $A^{[z],[\phi]} \to A^{[\phi]}$ and $A^{[z],[\phi]} \to A^{[z]}$, respectively. Then an isomorphism of extensions $\zeta \colon \mathcal{E}_{[\phi]}^{0,[z]} \to \mathcal{E}_{[z]}^{0,[\phi]}$ determines an isomorphism of extensions $\xi \colon \mathcal{E}^{z}_{[\phi]} \to \mathcal{E}_{[z]}^{\phi}$. Multiplying $\zeta$ by a character of $A^{[z],[\phi]}$ multiplies $\xi$ by the same character.
\end{proposition}

\begin{proof}
    The proof of \cite[Proposition 8.2]{Kal26b} uses calculations that pass word-for-word to the equal characteristic setting. The key new ingredient is this proof's use of Proposition 8.1 loc. cit., which makes extensive use of the duality isomorphism \eqref{eq:TNmap}, pairing \eqref{eq:TNbis}, and their formal properties developed in \S \ref{sec:mainsec} above. There is no new \v{C}ech-computation needed to pass the aforementioned proof of Proposition 8.1 to our setting---replacing the class $(z^{-1},t_{a}) \in Z^{1}(u \to W, Z \to T \xrightarrow{1-a} T)$ loc. cit. with $(z^{-1},t_{a},h) \in Z^{1}(\mathcal{E}, Z \to T \xrightarrow{1-a} T)$ (we warn the reader that our $h$ is different from the ``$h$'' appearing in \cite[\S 8.2]{Kal26b})---it relies entirely on the properties established in Proposition \ref{prop:pp1}, Theorem \ref{thm:TN}, and Corollary \ref{cor:TN}.
\end{proof}

We obtain the following result:

\begin{theorem}\label{thm:scuspidal}
    Let $G$ be a quasi-split connected reductive group over a local function field $F$ of characteristic $p>0$ which is split over a tamely ramified extension of $F$ and such that $p$ does not divide the order of the absolute Weyl group of $G$, with fixed finite central subgroup $Z \subset G$. Then the methods of \cite{Kal22} extend to give a construction of the $L$-packet $\Pi_{\phi}$ together with a bijection $\Pi_{\phi} 
    \to \textnormal{Irr}(\pi_{0}(Z_{\widehat{G}}(\phi) \times_{\widehat{G}} \widehat{G/Z}))$ as in \cite[(4.8)]{Kal22} for every discrete $L$-parameter $\phi \colon W_{F} \to {^{L}G}$. 
\end{theorem}

All of the representation-theoretic constructions from \cite{Kal22} are insensitive to the characteristic of the non-archimedean local field $F$ (and \cite{Kal22} does not explicitly assume mixed characteristic until it begins to use rigid inner forms). For the convenience of the reader, we make a list summarizing the constructions requiring rigid inner forms in \cite{Kal22} and their replacements in the function field setting. Throughout we fix a finite central subgroup $Z$ of $G$ and denote by $\psi$ an inner twist of $G$. For a supercuspidal $L$-parameter $\phi$ for $G$ we denote by $S$ the torus over $F$ obtained from $\phi$ as in \cite[\S 4.1]{Kal22}.

\begin{itemize}
    \item{In \cite[\S 4.2]{Kal22}, the notions of rigid inner twists and their representations from \cite[\S 5.1]{kaletha16} are replaced by the $Z$-rigid inner twists $(\psi,z',\bar{h})$ and their representations from \cite[Definition 7.1 and Definition 7.9]{Dillery1}. It follows from \cite[Proposition 5.12]{Dillery1} that if $Z$ contains $Z_{G_{\text{der}}}$ then every inner form of $G$ can be enriched to such a rigid inner twist.}
    \item{One can construct admissible embeddings $S \to G$ as in \cite[\S 4.2, \S 4.4]{Kal22} using \cite[Lemma 6.2]{Dillery1} and, for admissible embeddings $j_{0} \colon S \to G$ and $j \colon S \to G'$ for a $Z$-rigid inner twist and embedding $(G',\psi,z',\bar{h},j)$ of $G$, construct (using strongly regular elements of $S(F)$ and \cite[\S 7.1]{Dillery1}) an element $\text{inv}(j_{0},j) \in H^{1}(\mathcal{E}, Z \to S)$, which by \cite[Proposition 7.7, Lemma 7.8]{Dillery1} gives a bijection between the latter cohomology set and isomorphism classes of tuples $(G',\psi,z',\bar{h},j)$, replacing \cite[(4.3)]{Kal22}.}
    \item{The groups $S_{\phi}^{+}$, $\widehat{\bar{S}}$, and $[\widehat{\bar{S}}]^{\Gamma,+}$ are defined identically in our setting, and one has a Tate--Nakayama duality isomorphism \cite[Corollary 7.11]{Dillery1} identifying characters of $\pi_{0}([\widehat{\bar{S}}]^{\Gamma,+})$ with $H^{1}(\mathcal{E}, Z \to S)$ which is functorial in $[Z \to S]$.}
    \item{The result \cite[Lemma E.1]{Kal22} used for relating classes in $H^{1}(u \to W, Z_{G} \to S)$ to non-singular Deligne--Lusztig packets holds in our setting, whose proof is an elementary \v{C}ech cocycle computation.}
    \item{In the setting of \cite[\S 4.5]{Kal22}, we take the pair $(z,h)$ as fixed at the beginning of this section, $\phi$ essentially of depth zero, $S$, $\theta$, $\phi_{S}$, and $j_{0}$ as in loc. cit., and apply the framework of this section to the parameter $\phi_{S}$ (in place of $\phi$), $T =S$, $A=N_{G}(S)/S$ (with $S$ embedded in $G$ via $j_{0}$), $\tilde{T} = N_{G}(j_{0}(S))$, and $\bar{\mathcal{T}} = N_{\widehat{\bar{G}}}(\widehat{\bar{S}})$ equipped with the Galois action as explained in \cite{Kal22}. The two exactness assumptions from above hold in this setting by \cite[Lemma 3.4.10]{Kal19} and \cite[Lemma 4.5.3]{Kal22}, whose proofs hold in equal characteristic. We obtain a rigid inner twist $G'$ of $G$ and embedding $j$ of $S$ into $G'$ with invariant (as in the second bullet point) $[z]$ and have (using \cite[Proposition 7.7]{Dillery1}) $\tilde{T}_{z} = N_{G'}(jS)$.}
\end{itemize}

We can now prove the above theorem:
\begin{proof}
   The construction proceeds exactly as in \cite{Kal22} after replacing every cohomology set $H^{1}(u \to W, Z \to ?)$ (where $?$ can be any connected reductive group over $F$, such as $?=G$ or $?=T$ for an $F$-torus $T$) with $H^{1}(\mathcal{E}, Z \to ?)$ as defined in \cite{Dillery1} and the results in \cite{kaletha16} concerning each $H^{1}(u \to W, Z \to ?)$ with their analogues from \cite{Dillery1}. 

   More specifically, we use the first two bullet points above to construct the $L$-packet $\Pi_{\phi}$ as in \cite[\S 4.2]{Kal22}. The first four bullet points allow one to, using the arguments of \cite[\S 4.4]{Kal22}, reduce the claimed parametrization of $\Pi_{\phi}$ to the essentially depth zero case. In this case, after using \cite[Corollary 4.5.5, Corollary 4.5.7]{Kal22} along with the fifth bullet point we obtain the same key finite group extension comparison results \cite[Proposition 4.5.1, Lemma 4.5.8]{Kal22} after replacing the key input \cite[Proposition 8.2]{Kal26b} with our Proposition \ref{prop:prop8.2}. This proves the parametrization.
\end{proof}

\begin{remark}
    We defer any claims about the stability and transfer properties of the above parametrization of $\Pi_{\phi}$ as discussed in \cite[\S 4.6]{Kal22} to a later paper.
\end{remark}

\bibliography{FFdiscTNbib}
\bibliographystyle{amsalpha}

\end{document}